\documentclass[11pt,a4paper]{article}

\usepackage{latexsym}
\usepackage{mathrsfs}
\usepackage{graphicx}
\usepackage{supertabular}
\usepackage{subfigure}
\usepackage{setspace}
\usepackage{epstopdf}
\usepackage{amsmath,amsfonts,amsthm,amssymb}
\usepackage{float}
\usepackage{caption}
\usepackage[T1]{fontenc}
\usepackage{lmodern}
\usepackage{microtype}
\usepackage{mathtools}
\usepackage{enumitem}
\usepackage{xcolor}
\usepackage[colorlinks=true,linkcolor=blue!55!black,citecolor=blue!55!black,urlcolor=blue!55!black]{hyperref}
\usepackage[nameinlink,noabbrev]{cleveref}

\allowdisplaybreaks
\setlist{nosep}

\newtheorem{theorem}{Theorem}[section]
\newtheorem{lemma}[theorem]{Lemma}
\newtheorem{proposition}[theorem]{Proposition}
\newtheorem{corollary}[theorem]{Corollary}
\theoremstyle{definition}

\theoremstyle{remark}
\newtheorem{remark}[theorem]{Remark}

\theoremstyle{definition}

\newtheorem{conjecture}[theorem]{Conjecture}
\newtheorem{construction}[theorem]{Construction}
\newtheorem{problem}[theorem]{Problem}

\newtheorem{fact}[theorem]{Fact}

\newcommand{\Bn}{\mathcal{B}_n}
\newcommand{\cA}{\mathcal{A}}
\newcommand{\cF}{\mathcal{F}}
\newcommand{\cU}{\mathcal{U}}
\newcommand{\K}{\mathbb{K}}
\newcommand{\Q}{\mathbb{Q}}
\newcommand{\1}{\mathbf{1}}
\newcommand{\rank}{\operatorname{rank}}
\newcommand{\diag}{\operatorname{diag}}

\newcommand{\ol}[1]{\overline{#1}}

\newcommand{\binomcl}[1]{\binom{n}{\lfloor n/2\rfloor}}
\newcommand{\binompurdy}[1]{\binom{n}{\lfloor n/2\rfloor+1}}

\title{Hilton's Conjecture on Uncomplemented Incomparable Families}

\author{Yuejian Peng \footnote{School of Mathematics, Hunan University, Changsha 410082, P. R. China. E-mail: ypeng1@hnu.edu.cn.
 }~~~~ Zilong Yan \footnote{School of Mathematics, Hunan University,
 Changsha 410082, P. R. China. E-mail: zilongyan@hnu.edu.cn.} ~~~~ Cheng Zhao \footnote{SDU-ANU Joint Science College, Shandong University, Weihai, 264209, P. R. China. E-mail: chengzhao@sdu.edu.cn.} \\
 }

\date{October 1, 2026}

\begin{document}
\maketitle

\begin{abstract}
A family ${\mathcal A}\subseteq 2^{[n]}$ is \emph{uncomplemented} if it contains no set together with its complement. Families ${\mathcal A}_1,\ldots,{\mathcal A}_t$ are \emph{strongly incomparable} if every set in one family is incomparable with every set in another family.

In 1976, Hilton conjectured that if ${\mathcal A}_1,\ldots,{\mathcal A}_t\subseteq 2^{[n]}$ are uncomplemented and strongly incomparable, then
$$\sum_{i=1}^t |{\mathcal A}_i|\le 2^{n-1}.$$
We confirm this conjecture by a matrix rank argument that encodes incomparability and uncomplementedness in the same matrix, which is quite different from the earlier methods on this problem.

We also study weak incomparability, in which the same set may occur in different families. Goldwasser, Hilton and Zheng conjectured an upper bound in this setting and claimed that their conjecture is equivalent to Hilton's conjecture. We give counterexamples to their conjecture for every $t\ge3$ and identify the error in the claimed equivalence. We then propose a replacement conjecture and confirm it when $t$ is large.

The upper bound $2^{n-1}$ in Hilton's conjecture can be attained if some families are allowed to by empty.   It is natural to consider the maximum size of $\sum_{i=1}^t |{\mathcal A}_i|$ for uncomplemented and strongly incomparable families  ${\mathcal A}_1,\ldots,{\mathcal A}_t\subseteq 2^{[n]}$ when all families are required to be non-empty.  Denote this maximum number by $M(n,t)$.  We show that
\[
2^{n-1}-1\le M(n,t)\le 2^{n-1}
\]
for $n\ge3$ and $3\le t\le n$.


\end{abstract}

\section{Introduction and Background}\label{sec:introduction}

Let $[n]=\{1,2,\dots,n\}$ and let $\mathcal B_n=2^{[n]}$ denote the Boolean lattice of all subsets of $[n]$. 
For $X\subseteq[n]$ and ${\mathcal A}\subseteq 2^{[n]}$, define
\[
\overline X=[n]\setminus X,\qquad \overline{\mathcal A}=\{\overline X:X\in{\mathcal A}\}.
\]
A family ${\mathcal A}\subseteq 2^{[n]}$ is \emph{uncomplemented} if $X\in {\mathcal A}$ always implies $\overline X\notin {\mathcal A}$. Note that an uncomplemented family is a generalization of an intersecting family. Two sets $X,Y\subseteq[n]$ are \emph{strongly incomparable} if $X\ne Y$ and neither contains the other. They are \emph{weakly incomparable} if neither properly contains the other; in particular, equality is allowed. Families ${\mathcal A}_1,\dots,{\mathcal A}_t$ are called \emph{strongly}, respectively \emph{weakly}, \emph{incomparable} if every member of ${\mathcal A}_i$ is \emph{strongly}, respectively \emph{weakly}, \emph{incomparable} with every member of ${\mathcal A}_j$ whenever $i\ne j$.

A family ${\mathcal F}\subseteq 2^{[n]}$ is a \emph{Sperner family} (or an \emph{antichain}) if every  two distinct members of ${\mathcal F}$ are strongly incomparable.  Let us  recall several  classical results.

\begin{theorem}[Sperner--Purdy--Bollob\'as]\label{thm:classical-sperner-bounds}
Let ${\mathcal F}\subseteq2^{[n]}$ be a Sperner family. Then
\begin{enumerate}[label=\textup{(\roman*)}]
\item (Sperner \cite{Sperner1928}) $|{\mathcal F}|\le \binom{n}{\lfloor n/2\rfloor}$ ;
\item (Purdy \cite{Purdy1977}) if ${\mathcal F}$ is also uncomplemented, then
\[
|{\mathcal F}|\le \binom{n}{\lfloor n/2\rfloor+1};
\]
\item(Bollob\'as \cite{Bollobas1973}) if ${\mathcal F}$ is complemented, i.e. $F\in{\mathcal F}$ implies $\overline F\in{\mathcal F}$, then
\[
|{\mathcal F}|\le 2\binom{n-1}{\lfloor n/2\rfloor-1}.
\]
\end{enumerate}
\end{theorem}
 Seymour proved the following result for two  non-empty, strongly incomparable families.

\begin{theorem}[Seymour \cite{Seymour1973}]\label{thm:seymour-background}
If ${\mathcal A}$ and ${\mathcal B}$ are two non-empty, strongly incomparable families of subsets of $[n]$, then
\[
\sqrt{|{\mathcal A}|}+\sqrt{|{\mathcal B}|}\le 2^{n/2}.
\]
\end{theorem}

For several weakly incomparable families, Daykin, Frankl, Greene and Hilton obtained a direct generalization of Sperner's theorem.

\begin{theorem}[Daykin--Frankl--Greene--Hilton \cite{DFGH1981}]\label{thm:dfgh-background}
Let ${\mathcal A}_1,\ldots,{\mathcal A}_t$ be mutually weakly incomparable families of subsets of $[n]$. Then
\[
\sum_{i=1}^t|{\mathcal A}_i|
\le \max\left\{2^n,\;t\binom{n}{\lfloor n/2\rfloor}\right\}.
\]
\end{theorem}


The classical extremal problem for uncomplemented families begins with the elementary observation that a single uncomplemented family has size at most $2^{n-1}$, since the Boolean lattice is partitioned into complementary pairs. Hilton proved that the same bound remains valid for two families, provided they are also strongly incomparable.

\begin{theorem}[Hilton \cite{Hilton1976}, 1976] \label{thm:Hilton-two}
If ${\cal A}_{1}, {\cal A}_{2}\subseteq 2^{[n]}$ are uncomplemented and strongly incomparable,
 then
\[
|{\cal A}_1|+|{\cal A}_2| \le 2^{n-1}.
\]
\end{theorem}
The bound in Hilton's theorem is sharp, as shown by the following construction.

\begin{construction}\label{cons-h}
Let $1 \le i \le n-1$ and set
 $$ {\cal A}_1 = \{X \subseteq [n] : 1 \in X, |X| \le i\}, $$
 $$ {\cal A}_2 = \{X \subseteq [n] : 1 \notin X, |X| \ge i\}.$$
\end{construction}

In his \emph{Mathematical Reviews} review of Hilton's paper, Kleitman \cite{KleitmanMR53-146} gave an alternative proof of Theorem~\ref{thm:Hilton-two} using a correlation inequality and also emphasized the following conjecture, which Hilton had proposed in 1976.

\begin{conjecture}[Hilton's Conjecture \cite{KleitmanMR53-146}] \label{ch}
If ${\cal A}_{1}, {\cal A}_{2}, \dots, {\cal A}_t\subseteq 2^{[n]}$ are uncomplemented and strongly incomparable, then
\[
\sum_{i=1}^t|{\cal A}_i| \le 2^{n-1}.
\]
\end{conjecture}

This beautiful conjecture combines complementary property (a generalization of intersecting property) and incomparable property (related to the concept of antichain in classical Sperner's theorem). Liu--Zhao \cite{LiuZhao2001} established an upper bound for $\sum_{i=1}^t|{\cal A}_i|$. Peng--Zhao \cite{PengZhao2005} subsequently improved this bound to
$\sum_{i=1}^t|{\cal A}_i| < (4-2\sqrt{2})2^{n-1}$. The general conjecture has  remained open for several decades.
In this paper, we confirm Hilton's  conjecture, i.e, prove the following result.

\begin{theorem}\label{thm:hilton}
If ${\mathcal A}_1,\dots,{\mathcal A}_t\subseteq 2^{[n]}$ are uncomplemented and strongly incomparable, then
\[
\sum_{i=1}^t |{\mathcal A}_i|\le 2^{n-1}.
\]
\end{theorem}
The bound $2^{n-1}$ is attained by taking ${\mathcal A}_1$ and ${\mathcal A}_2$ as in Construction~\ref{cons-h} and taking all remaining families to be empty.

The main new ingredient in the proof is a matrix rank argument that encodes incomparability and uncomplementedness simultaneously. We construct a matrix $M$, indexed by the Boolean lattice, with rank $2^{n-1}$. Its entries satisfy $M_{A,B}=0$ whenever $A$ and $B$ are incomparable, while $M[\mathcal{F},\mathcal{F}]$ is nonsingular for every uncomplemented family $\mathcal{F}$. Thus, for strongly incomparable families, the relevant principal submatrix becomes block diagonal, and Hilton's conjecture follows from a direct rank comparison. This method is different from the correlation and probabilistic methods used in earlier work.

The situation changes when strong incomparability is replaced by weak incomparability. In this setting, the same set may appear in several families, and this new possibility creates additional extremal configurations. Goldwasser, Hilton and Zheng \cite{GHZ19} proposed the following bound.

\begin{conjecture}[ Goldwasser, Hilton and Zheng \cite{GHZ19}] \label{GHZ}
If ${\mathcal A}_1, \dots, {\mathcal A}_t\subseteq 2^{[n]}$ are uncomplemented and weakly incomparable, then
\[
|{\mathcal A}_1|+\cdots+|{\mathcal A}_t|\le\max\left\{2^{n-1}, t\binom{n}{\lfloor\frac{n}{2}\rfloor+1}\right\}.
\]
\end{conjecture}

Our second contribution concerns weak incomparability. We show that the conjecture of Goldwasser, Hilton and Zheng is false for every $t\ge3$. The counterexamples also show that their claimed equivalence with Hilton's conjecture cannot hold. In Appendix~\ref{app:ghz-error}, we identify the precise step in their proof where the argument fails. Motivated by two natural constructions, we propose a replacement conjecture and prove it in the range
\[
t\binom{n}{\lfloor n/2\rfloor+1}\ge 2^n-2.
\]
The proof uses a reduction for repeated sets together with a weighted LYM inequality.

We now give two constructions showing that Conjecture~\ref{GHZ} is not correct once $t\ge3$. For a set $A$ and a non-negative integer $k$, let ${A \choose k}$ denote the family of all $k$-subsets of $A$.

\smallskip\noindent
\emph{Construction A.}
Let $1\le h\le n-1$. Define
\begin{eqnarray*}
{\mathcal A}_1&=&\{X\subseteq[n]: 1\in X,\ |X|\le n-h\},\\
{\mathcal A}_2&=&\cdots={\mathcal A}_t=\binom{[n]\setminus\{1\}}{n-h}.
\end{eqnarray*}
These families are uncomplemented and weakly incomparable. Also,
\[
|{\mathcal A}_1|+\cdots+|{\mathcal A}_t|
=2^{n-1}-\sum_{i=0}^{h-1}\binom{n-1}{i}+(t-1)\binom{n-1}{h-1}.
\]
In particular, taking $h=1$ gives examples with total size $2^{n-1}+t-2$.

\smallskip\noindent
\emph{Construction B.}
Take
\[
{\mathcal T}=\binom{[n]}{\lfloor n/2\rfloor+1}.
\]
This is an uncomplemented antichain. Set ${\mathcal A}_1=\cdots={\mathcal A}_t={\mathcal T}$. Then
\[
|{\mathcal A}_1|+\cdots+|{\mathcal A}_t|
=t\binom{n}{\lfloor n/2\rfloor+1}.
\]

These two constructions give the following conclusion.

\begin{fact}\label{thm:mainweaklower}
There exist uncomplemented and weakly incomparable families ${\mathcal A}_1,\dots,{\mathcal A}_t\subseteq 2^{[n]}$ such that
\[
\sum_{i=1}^t |{\mathcal A}_i|\ge
\max_{1\le h\le n-1}\left\{2^{n-1}-\sum_{i=0}^{h-1}\binom{n-1}{i}+(t-1)\binom{n-1}{h-1},\,
 t\binom{n}{\lfloor\frac{n}{2}\rfloor+1}\right\}.
\]
\end{fact}

\begin{remark}
Taking $h=1$ in Fact~\ref{thm:mainweaklower}  gives
\[
\sum_{i=1}^t |{\mathcal A}_i|\ge
\max\left\{2^{n-1}+(t-2),  t\binom{n}{\lfloor\frac{n}{2}\rfloor+1}\right\}.
\]
For each fixed $t\ge3$, choose $n$ sufficiently large that
$t\binom{n}{\lfloor n/2\rfloor+1}\le2^{n-1}$; this is possible since
$\binom{n}{\lfloor n/2\rfloor+1}=O(2^n/\sqrt n)$.
The first term above then exceeds the bound in Conjecture~\ref{GHZ}.
Thus that conjecture fails for every $t\ge3$.
\end{remark}

 Constructions $A$ and $B$ lead us to the following replacement for Conjecture~\ref{GHZ}.

\begin{conjecture}\label{conj:weak}
Let $n\ge2$ and $t\ge2$. If ${\mathcal A}_1,\dots,{\mathcal A}_t\subseteq 2^{[n]}$ are uncomplemented and weakly incomparable, then
\[
\sum_{i=1}^t |{\mathcal A}_i| \le
\max_{1\le h\le n-1}\left\{2^{n-1}-\sum_{i=0}^{h-1}\binom{n-1}{i}+(t-1)\binom{n-1}{h-1},\,
 t\binom{n}{\lfloor\frac{n}{2}\rfloor+1}\right\}.
\]
\end{conjecture}

We show the above conjecture holds when  $t$ is sufficiently large.

\begin{theorem}\label{partial-conj}
Let $n\ge 2$ and $t\ge2$ be integers such that $t\binom{n}{\lfloor n/2\rfloor+1}\ge2^n-2$.
If ${\mathcal A}_1,\ldots,{\mathcal A}_t\subseteq2^{[n]}$ are
uncomplemented and weakly incomparable, then
\[
  \sum_{i=1}^t|{\mathcal A}_i|\le t\binom{n}{\lfloor n/2\rfloor+1}.
\]
Equality is attained by taking
\[
  {\mathcal A}_1=\cdots={\mathcal A}_t
  =\binom{[n]}{\lfloor n/2\rfloor+1}.
\]
Consequently, Conjecture~\ref{conj:weak} holds whenever $n\ge 2$ and
$t\binom{n}{\lfloor n/2\rfloor+1}\ge2^n-2$.
\end{theorem}

Note that the bound $2^{n-1}$ in Theorem \ref{thm:hilton}  is attained by taking ${\mathcal A}_1$ and ${\mathcal A}_2$ as in Construction~\ref{cons-h} and taking all remaining families to be empty. Is the bound $2^{n-1}$ still reachable if we require that all families are non-empty?
   Let $M(n,t)$
be the maximum of $\sum_{i=1}^t|{\mathcal A}_i|$ when all $t$ families
are non-empty, uncomplemented and strongly incomparable. It is easy to see that
$M(n,1)=2^{n-1}$, and  $M(n,2)=2^{n-1}$ by Hilton's theorem  for $n\ge2$. We show  that
\[
  2^{n-1}-1\le M(n,t)\le2^{n-1}
  \qquad(n\ge3,\ 3\le t\le n).
\]
We conjecture that $M(n,t)=2^{n-1}-1$ if $n\ge3$ and $3\le t\le n$.

We prove  Theorem~\ref{thm:hilton}  in Section~\ref{sec2}.
In Section~\ref{sec3}, we  prove
Theorem~\ref{partial-conj} by using a reduction of repeated sets and a
weighted LYM inequality. In Section~\ref{sec:concluding}, we discuss the
non-empty version of Hilton's problem and states a related conjecture.
Appendix~\ref{app:ghz-error} explains the error in the proof of
Theorem~1.5 of Goldwasser--Hilton--Zheng \cite{GHZ19}.

\section{Proof of Hilton's Conjecture on Strong Incomparability}\label{sec2}

In this section, we will prove Theorem~\ref{thm:hilton}. Recall that $\mathcal B_n=2^{[n]}$ denote the Boolean lattice of all subsets of $[n]$.
Let $\Omega_n$ denote the set of complementary orbits
\[
  \Omega_n=\bigl\{\{A,\ol A\}:A\in\Bn\bigr\}.
\]
Since $n\ge1$, every orbit has two elements and $|\Omega_n|=2^{n-1}$.

For each $\omega\in\Omega_n$, choose one of its two elements and denote it by $A_\omega$. Let $t_\omega$ be an indeterminate, with the indeterminates $\{t_\omega:\omega\in\Omega_n\}$ algebraically independent over $\Q$. Define
\[
  \K=\Q(t_\omega^{\pm1}:\omega\in\Omega_n)
\]
and, for every $A\in\Bn$, define $q_A\in\K$ by
\[
  q_{A_\omega}=t_\omega,
  \qquad
  q_{\ol{A_\omega}}=t_\omega^{-1}.
\]
Thus
\begin{equation}\label{eq:reciprocal}
  q_Aq_{\ol A}=1
  \qquad (A\in\Bn).
\end{equation}

Index all matrices below by $\Bn$. Define the matrix $Z$ and the  matrix $P$ as follows.
\begin{equation}\label{eq:ZP}
  Z_{A,B}=\1_{\{A\subseteq B\}},
  \qquad
  P_{A,B}=\1_{\{A=\ol B\}}.
\end{equation}
Let $T$ be the $2^n\times 2^n$ matrix indexed by $\Bn$ and defined by
\begin{equation}\label{eq:inverse}
  T_{A,B}
  =(-1)^{|B\setminus A|}\1_{\{A\subseteq B\}}.
\end{equation}
Then
\begin{align}
  (ZT)_{A,B}
  &=\sum_{C\in \Bn} Z_{A,C}T_{C,B} \notag\\
  &=\sum_{C\in \Bn}
    \1_{\{A\subseteq C\}}
    (-1)^{|B\setminus C|}
    \1_{\{C\subseteq B\}} \notag\\
  &=\sum_{A\subseteq C\subseteq B}(-1)^{|B\setminus C|}
   =\1_{\{A=B\}}.\label{eq:times}
\end{align}
Indeed, the two indicator factors in the preceding sum force precisely
$A\subseteq C\subseteq B$. If $A\nsubseteq B$, there is no such $C$, so the
sum is zero. If $A\subseteq B$, putting $S=B\setminus C$ gives a bijection
between the sets $C$ with $A\subseteq C\subseteq B$ and the subsets
$S\subseteq B\setminus A$. Therefore
\[
  \sum_{A\subseteq C\subseteq B}(-1)^{|B\setminus C|}
  =\sum_{S\subseteq B\setminus A}(-1)^{|S|}
  =(1-1)^{|B\setminus A|},
\]
which is $1$ when $A=B$ and $0$ otherwise. Thus $ZT=I$, and hence
\begin{equation}\label{eq:mobius}
  (Z^{-1})_{A,B}
  =(-1)^{|B\setminus A|}\1_{\{A\subseteq B\}}.
\end{equation}
Let
\begin{equation}\label{eq:defM}
  D=\diag(q_A:A\in\Bn),
  \qquad
  M=Z^{-1}(D+P)Z.
\end{equation}

We will prove the following theorem and apply it to prove Theorem \ref{thm:hilton}.

\begin{theorem}\label{thm:certificate}
The matrix $M$ in \eqref{eq:defM} has the following properties.
\begin{enumerate}[label=\textup{(\roman*)}]
\item $\rank_{\K}M=2^{n-1}$.
\item If $A,B\in\Bn$ are incomparable, then $M_{A,B}=0$.
\item If $\cF\subseteq\Bn$ is uncomplemented, then $M[\cF,\cF]$ is nonsingular over $\K$.
\end{enumerate}
\end{theorem}

We prove the three assertions in the next  section.

\subsection{Proof of Theorem \ref{thm:certificate}}

We first show Part (i) of Theorem \ref{thm:certificate}.

\begin{proposition}\label[proposition]{prop:rank}
The matrix $M$ has rank $2^{n-1}$.
\end{proposition}

\begin{proof}
Order the rows and columns of $D+P$ orbit by orbit. On the ordered orbit $(A_\omega,\ol{A_\omega})$, the corresponding block is
\[
  R_\omega=\begin{pmatrix}
    t_\omega & 1\\
    1 & t_\omega^{-1}
  \end{pmatrix}
  =
  \begin{pmatrix}t_\omega\\1\end{pmatrix}
  \begin{pmatrix}1&t_\omega^{-1}\end{pmatrix}.
\]
Each block has rank one. Thus
\begin{equation}\label{eq:orbit-blocks}
  D+P=\diag(R_\omega:\omega\in\Omega_n).
\end{equation}
There are $2^{n-1}$ complementary orbits, so
\[
  \rank(D+P)=2^{n-1}.
\]
Since $Z$ is invertible, conjugation preserves rank, and hence
\[
  \rank M=2^{n-1}.
\]
\end{proof}

We next determine the zero pattern of $M$. Write
\begin{equation}\label{eq:UJ}
  M=U+J,
  \qquad
  U=Z^{-1}DZ,
  \qquad
  J=Z^{-1}PZ.
\end{equation}

\begin{lemma}\label[lemma]{lem:U}
For $A,B\in\Bn$,
\begin{equation}\label{eq:Uentry}
  U_{A,B}
  =
  \begin{cases}
  \displaystyle
  \sum_{A\subseteq C\subseteq B}(-1)^{|C\setminus A|}q_C,
      & A\subseteq B,\\[3mm]
  0,  & A\nsubseteq B.
  \end{cases}
\end{equation}
In particular, $U_{A,B}$ can be nonzero only if $A\subseteq B$.
\end{lemma}

\begin{proof}
Using \eqref{eq:ZP} and \eqref{eq:mobius},
\[
  U_{A,B}
  =\sum_{C\in\Bn}(Z^{-1})_{A,C}q_C Z_{C,B}
  =\sum_{A\subseteq C\subseteq B}(-1)^{|C\setminus A|}q_C.
\]
The indexing interval is empty unless $A\subseteq B$.
\end{proof}

\begin{lemma}\label[lemma]{lem:J}
For $A,B\in\Bn$,
\begin{equation}\label{eq:Jentry}
  J_{A,B}=(-1)^{n-|A|}\1_{\{B\subseteq A\}}.
\end{equation}
In particular, $J_{A,B}$ can be nonzero only if $B\subseteq A$.
\end{lemma}

\begin{proof}
From \eqref{eq:ZP} and \eqref{eq:mobius},
\begin{align*}
  J_{A,B}
  &=\sum_{C,E\in\Bn}(Z^{-1})_{A,C}P_{C,E}Z_{E,B}\\
  &=\sum_{\substack{C\supseteq A\\ \ol C\subseteq B}}
      (-1)^{|C\setminus A|}\\
  &=\sum_{C\supseteq A\cup\ol B}(-1)^{|C\setminus A|}.
\end{align*}
Put $S=A\cup\ol B$. Since $A\subseteq S$, every $C\supseteq S$ is uniquely of the form $C=S\cup T$ with $T\subseteq[n]\setminus S$, and hence
\begin{align*}
  \sum_{C\supseteq S}(-1)^{|C\setminus A|}
  &=(-1)^{|S\setminus A|}
    \sum_{T\subseteq[n]\setminus S}(-1)^{|T|}\\
  &=(-1)^{|S\setminus A|}(1-1)^{n-|S|}.
\end{align*}
Thus the sum is zero unless $S=[n]$. The condition $S=[n]$ is equivalent to $B\subseteq A$; in that case the sum has the single term $C=[n]$ and equals $(-1)^{n-|A|}$.
\end{proof}

Now we can get Part (ii) of Theorem \ref{thm:certificate}.

\begin{corollary}\label[corollary]{cor:support}
If $A$ and $B$ are incomparable, then $M_{A,B}=0$.
\end{corollary}

\begin{proof}
If $A$ and $B$ are incomparable, then neither $A\subseteq B$ nor $B\subseteq A$. Apply Lemmas \ref{lem:U} and \ref{lem:J} to \eqref{eq:UJ}.
\end{proof}

Note that the entries of $M$ may be written explicitly as
\begin{equation}\label{eq:Mexplicit}
M_{A,B}=
\begin{cases}
\displaystyle
\sum_{A\subseteq C\subseteq B}(-1)^{|C\setminus A|}q_C,
  & A\subsetneq B,\\[3mm]
q_A+(-1)^{n-|A|},
  & A=B,\\[1mm]
(-1)^{n-|A|},
  & B\subsetneq A,\\[1mm]
0,
  & A\text{ and }B\text{ are incomparable}.
\end{cases}
\end{equation}

Next we show Part (iii) of Theorem \ref{thm:certificate}.
 Let
\[
  \mathbf y=(y_A:A\in\Bn)
\]
be algebraically independent commuting variables, and define the
$\Q$-algebra homomorphism
\begin{equation}\label{eq:phi}
  \phi:\Q[\mathbf y]\longrightarrow
  \Q[t_\omega^{\pm1}:\omega\in\Omega_n],
  \qquad
  \phi(y_A)=q_A.
\end{equation}
For a monomial $\vec t^{\  \vec a}=\prod_\omega t_\omega^{a_\omega}$, define its
\emph{orbit weight} by
\[
  \operatorname{wt}(\vec t^{\  \vec a})=\|\vec a\|_1
  =\sum_{\omega\in\Omega_n}|a_\omega|.
\]

\begin{lemma}\label[lemma]{lem:weight-separation}
Let $\cF\subseteq\Bn$ be uncomplemented, let $m=|\cF|$, and suppose that
\[
  p(\mathbf y)
  =\prod_{A\in\cF}y_A+r(\mathbf y),
  \qquad
  \deg r\le m-1.
\]
Then $\phi(p)\ne0$.
\end{lemma}

\begin{proof}
Because $\cF$ is uncomplemented, it meets each complementary orbit in at
most one set. Consequently
\[
  \phi\left(\prod_{A\in\cF}y_A\right)
  =\prod_{A\in\cF}q_A
\]
is a monomial involving exactly $m$ distinct orbit variables,
each with exponent $1$ or $-1$. Its orbit weight is therefore $m$.

Let $\mu=\prod_A y_A^{\alpha_A}$ be a monomial occurring in
$r$, and put $d=\deg\mu=\sum_A\alpha_A\le m-1$. For each
$A\in\Bn$, write
\[
  q_A=t_{\omega(A)}^{\varepsilon_A},
  \qquad \varepsilon_A\in\{1,-1\},
\]
where $\omega(A)=\{A,\ol A\}$. If
\[
  \phi(\mu)=\prod_{\omega\in\Omega_n}t_\omega^{b_\omega},
\]
then
\[
  b_\omega=\sum_{A\in\omega}\varepsilon_A\alpha_A.
\]
Consequently, by the triangle inequality,
\[
  \operatorname{wt}(\phi(\mu))
  =\sum_\omega|b_\omega|
  \le\sum_\omega\sum_{A\in\omega}\alpha_A
  =d
  \le m-1.
\]
Thus no monomial of $\phi(r)$ can coincide with
$\prod_{A\in\cF}q_A$. The latter has coefficient $1$ in $\phi(p)$ and
cannot be cancelled, so $\phi(p)\ne0$.
\end{proof}

\begin{proposition}\label[proposition]{prop:minor}
For every uncomplemented family $\cF\subseteq\Bn$,
\[
  \det M[\cF,\cF]\ne0.
\]
\end{proposition}

\begin{proof}
The assertion is immediate when $\cF=\varnothing$, so assume that
$m=|\cF|\ge1$. Put
\[
  Y=\diag(y_A:A\in\Bn),
  \qquad
  \widetilde M=Z^{-1}(Y+P)Z,
\]
and write
\[
  \widetilde U=Z^{-1}YZ,
  \qquad
  J=Z^{-1}PZ,
\]
so that $\widetilde M=\widetilde U+J$. Define
\[
  p_{\cF}(\mathbf y)
  =\det \widetilde M[\cF,\cF]
  \in\Q[\mathbf y].
\]
Every entry of $\widetilde U$ is homogeneous of degree one in the
variables $y_A$, whereas $J$ is independent of them. Consequently
$p_{\cF}$ has total degree at most $m$, and its homogeneous component of
degree $m$ is
\begin{equation}\label{eq:topcomponent}
  \det \widetilde U[\cF,\cF].
\end{equation}

Choose a linear extension of the inclusion order on $\Bn$, so that $A\subseteq B$ implies that $A$ appears no later than $B$. Since
\[
  Z_{A,B}=\1_{\{A\subseteq B\}},
\]
all nonzero entries of $Z$ lie on or above the diagonal, and hence $Z$ is upper triangular. By the M\"obius formula \eqref{eq:mobius}, the same is true for $Z^{-1}$. Since $Y$ is diagonal, it follows that $\widetilde U=Z^{-1}YZ$ is also upper triangular. Moreover,
\[
  \widetilde U_{A,A}
  =\sum_{C\in\Bn}(Z^{-1})_{A,C}y_CZ_{C,A}
  =y_A,
\]
because a nonzero summand requires $A\subseteq C\subseteq A$, and therefore $C=A$.
The principal submatrix $\widetilde U[\cF,\cF]$, with the induced
ordering, is therefore upper triangular with diagonal
$(y_A:A\in\cF)$. Hence
\begin{equation}\label{eq:unique-top}
  \det \widetilde U[\cF,\cF]
  =\prod_{A\in\cF}y_A.
\end{equation}
It follows from \eqref{eq:topcomponent} and \eqref{eq:unique-top} that
\begin{equation}\label{eq:pdecomp}
  p_{\cF}(\mathbf y)
  =\prod_{A\in\cF}y_A+r_{\cF}(\mathbf y),
  \qquad
  \deg r_{\cF}\le m-1.
\end{equation}

Applying Lemma \ref{lem:weight-separation} to \eqref{eq:pdecomp} gives $\phi(p_{\cF})\ne0$. Since applying
$\phi$ entrywise to $\widetilde M$ gives $M$, and the determinant is a
polynomial in the matrix entries,
\[
  \phi(p_{\cF})=\det M[\cF,\cF].
\]
Therefore $\det M[\cF,\cF]\ne0$.
\end{proof}

Combining Proposition \ref{prop:rank}, Corollary \ref{cor:support}, and
Proposition \ref{prop:minor}, we obtain Theorem \ref{thm:certificate}.

\subsection{Proof of  Theorem \ref{thm:hilton}}

In this section, we will apply Theorem \ref{thm:certificate} to prove Theorem \ref{thm:hilton}.

\begin{proof}[Proof of Theorem \ref{thm:hilton}]
Let
\[
  \cU=\bigcup_{i=1}^t\cA_i.
\]
Strong incomparability forbids equality between sets belonging to distinct families, so this union is disjoint and
\begin{equation}\label{eq:disjoint}
  |\cU|=\sum_{i=1}^t|\cA_i|.
\end{equation}

Fix the matrix $M$ from  Theorem \ref{thm:certificate}. If $i\ne j$, $A\in\cA_i$, and $B\in\cA_j$, then $A$ and $B$ are incomparable. By Theorem \ref{thm:certificate}(ii),
\[
  M_{A,B}=0.
\]
Consequently, after ordering the indices family by family,
\begin{equation}\label{eq:block}
  M[\cU,\cU]
  =\diag\bigl(
       M[\cA_1,\cA_1],\ldots,
       M[\cA_t,\cA_t]
    \bigr).
\end{equation}
Each $\cA_i$ is uncomplemented, so Theorem \ref{thm:certificate}(iii) shows that every diagonal block in \eqref{eq:block} is nonsingular. Hence
\[
  \rank M[\cU,\cU]
  =\sum_{i=1}^t|\cA_i|.
\]
A submatrix cannot have larger rank than the ambient matrix. By Theorem \ref{thm:certificate}(i),
\[
  \sum_{i=1}^t|\cA_i|
  =\rank M[\cU,\cU]
  \le\rank M
  =2^{n-1}.
\]
\end{proof}

\section{Results on Weak Incomparability}\label{sec3}

In this section, we prove Theorem~\ref{partial-conj}. We first show the following property.

\begin{lemma}\label{lem:weak-saturation}
Let $n\ge1$ and $t\ge2$, and suppose that
${\mathcal A}_1,\ldots,{\mathcal A}_t\subseteq2^{[n]}$ are
uncomplemented and weakly incomparable. There exist uncomplemented,
weakly incomparable families
\[
  {\mathcal A}'_i={\mathcal T}\mathbin{\dot\cup}{\mathcal U}_i
  \qquad(1\le i\le t)
\]
with
\[
  \sum_{i=1}^t|{\mathcal A}'_i|
  \ge\sum_{i=1}^t|{\mathcal A}_i|,
\]
such that the following hold:
\begin{enumerate}[label=\textup{(\roman*)}]
\item ${\mathcal T}$ is an uncomplemented antichain.
\item The families ${\mathcal U}_1,\ldots,{\mathcal U}_t$ are pairwise
      disjoint, uncomplemented and strongly incomparable.
\item Writing ${\mathcal U}=\bigcup_{i=1}^t{\mathcal U}_i$, every member
      of ${\mathcal T}$ is strongly incomparable with every member of
      ${\mathcal U}$, and
      $\overline{\mathcal T}\cap{\mathcal U}=\varnothing$.
\end{enumerate}
\end{lemma}

\begin{proof}
Let
\[
  m(X)=|\{i\in[t]:X\in{\mathcal A}_i\}|.
\]
If $m(X)\ge2$, then $X$ is incomparable with every other set
in the system. Indeed, if $Y\ne X$ occurs in ${\mathcal A}_j$, we can take
$i\ne j$ with $X\in{\mathcal A}_i$, and weak incomparability implies that $X$ and $Y$ are incomparable.

Suppose $2\le m(X)<t$. Since each family is uncomplemented,
\[
  m(X)+m(\overline X)\le t.
\]
Delete $\overline X$ from every family and add $X$ to every family.
The observation above guarantees that weak incomparability is preserved.
Deleting $\overline X$ makes sure  that no complementary pair is introduced.
The change in total size is
\[
  t-m(X)-m(\overline X)\ge0.
\]
This operation makes the multiplicities on the complementary pair
$\{X,\overline X\}$ equal to $(t,0)$; later operations on other pairs
do not change them. Thus after finitely many operations every set in the system has
multiplicity  $1$ or $t$.

Let ${\mathcal T}$ consist of the sets of multiplicity $t$, and let
${\mathcal U}_i$ consist of the sets of multiplicity $1$ in the $i$-th
family. The initial observation shows that ${\mathcal T}$ is an
antichain and is strongly incomparable with ${\mathcal U}$.
The families ${\mathcal U}_i$ are disjoint, so their weak incomparability
is strong incomparability. Finally, ${\mathcal T}$ is contained in every
${\mathcal A}'_i$, each of which is uncomplemented. This proves both
that ${\mathcal T}$ is uncomplemented and that no member of
${\mathcal U}$ is the complement of a member of ${\mathcal T}$.
\end{proof}

\subsection{A weighted LYM inequality}
 For $n\ge1$, let
\begin{equation}\label{eq:weak-parameters}
  B_n=\binom{n}{\lfloor n/2\rfloor+1}.
\end{equation}
For $0\le j\le n$, define
\[
  d_n(j)=\begin{cases}
    B_n, & n\text{ is even and }j=n/2,\\
    \binom{n}{j}, & \text{otherwise},
  \end{cases}
  \qquad
  L_n({\mathcal F})=\sum_{F\in{\mathcal F}}\frac1{d_n(|F|)}.
\]
We first show the following.
\begin{lemma}
\label{lem:unpaired-lym}
Let $n=2k\ge2$, and let ${\mathcal F}\subseteq2^{[n]}$ be any antichain.
Define
\[
  P({\mathcal F})
  =\{F\in{\mathcal F}:|F|=k,\ \overline F\notin{\mathcal F}\}.
\]
Then
\begin{equation}\label{eq:unpaired-lym}
  \sum_{F\in{\mathcal F}}\frac1{\binom{n}{|F|}}
  +\frac{|P({\mathcal F})|}{k\binom{n}k}\le1.
\end{equation}
\end{lemma}

\begin{proof}
If ${\mathcal F}$ contains $\varnothing$ or $[n]$, it has only one
member and the conclusion is immediate. Otherwise, fix an oriented
cyclic ordering of $[n]$. Let $a$ count the members of ${\mathcal F}$
that are non-empty proper intervals, and let $b$ count those intervals
that belong to $P({\mathcal F})$.

An antichain of proper intervals has at most $n$ members, since their
starting points are distinct. Moreover, if $a=n$, then all these intervals have the same length. Indeed, if every starting point occurs, write
$\ell_i$ for the length of the interval starting at position $i$, with
indices taken modulo $n$. If $\ell_i>\ell_{i+1}$, the interval starting
at $i+1$ is properly contained in the one starting at $i$. Thus
$\ell_i\le\ell_{i+1}$ for every $i$, forcing all lengths to be equal.

The $n$ intervals of length $k$ form $k$ complementary pairs. At most
one interval from each pair belongs to $P({\mathcal F})$, so $b\le k$.
If $b>0$, equality $a=n$ would force all $n$ intervals of length $k$ to
belong to ${\mathcal F}$. Each of them would then have its complement
in ${\mathcal F}$, contradicting $b>0$. Hence $a\le n-1$ whenever
$b>0$, and in all cases
\[
  a+\frac bk\le n.
\]
A fixed $j$-set, $1\le j\le n-1$, is an interval in a uniformly random
cyclic ordering with probability $n/\binom{n}{j}$. Taking expectations
and dividing by $n$ proves \eqref{eq:unpaired-lym}.
\end{proof}

\begin{corollary}\label{lem:uncomplemented-lym}
If ${\mathcal F}\subseteq2^{[n]}$ is an uncomplemented antichain, then
\begin{equation}\label{eq:weighted-lym}
  L_n({\mathcal F})\le1.
\end{equation}
Consequently, $|{\mathcal F}|\le B_n$.
\end{corollary}

\begin{proof}
For odd $n$, this is the usual Lubell--Yamamoto--Meshalkin (LYM)
inequality \cite{Yamamoto1954,Meshalkin1963,Lubell1966}. For completeness,
we recall the standard chain interpretation. A maximal chain in the
Boolean lattice has the form
\[
  \varnothing=C_0\subset C_1\subset\cdots\subset C_n=[n],
  \qquad |C_j|=j.
\]
Choose a maximal chain uniformly at random; equivalently, choose a
uniformly random permutation $\pi$ of $[n]$ and put
$C_j=\{\pi(1),\ldots,\pi(j)\}$. If $F$ is a fixed $j$-set, then exactly
$j!(n-j)!$ of the $n!$ permutations generate a chain containing $F$ at
level $j$. Hence
\[
  \Pr(F\in\mathcal C)=\frac{j!(n-j)!}{n!}
  =\binom{n}{j}^{-1}.
\]
Since an antichain meets any maximal chain in at most one member,
linearity of expectation gives
\[
  \sum_{F\in\mathcal F}\frac1{\binom{n}{|F|}}
  =\mathbb E\,|\mathcal C\cap\mathcal F|\le1.
\]
This is the usual LYM inequality and yields the asserted odd-$n$ case.

For $n=2k$, every middle-level member of ${\mathcal F}$ belongs to
$P({\mathcal F})$. Thus each
middle-level member contributes
\[
\frac1{d_n(k)}=\frac1{B_n}=\frac1{\binom{2k}{k+1}}
 =\ \frac1{\binom{2k}k}+\frac1{k\binom{2k}k}
\]
to $L_n({\mathcal F})$,
while every member on another level $j$ contributes
$$1/\binom{2k}j=1/d_n(j)$$ to $L_n({\mathcal F})$. In view of (\ref{eq:unpaired-lym}) and   Lemma~\ref{lem:unpaired-lym}, we
have $L_n({\mathcal F})\le1$ when $n$ is even.

Finally, since $d_n(j)\le B_n$ for every $0\le j\le n$, we have
\[
  \frac{|{\mathcal F}|}{B_n}\le L_n({\mathcal F})\le1,
\]
and hence $|{\mathcal F}|\le B_n$.
\end{proof}

\subsection{Proof of Theorem~\ref{partial-conj}}\label{weakupper}

\begin{proof}[Proof of Theorem~\ref{partial-conj}]
Apply Lemma~\ref{lem:weak-saturation}. Since the resulting system has total
size at least as large as the original one, it suffices to prove the
bound for the resulting families. Retain the notation
${\mathcal T},{\mathcal U}_1,\ldots,{\mathcal U}_t$ from that lemma,
and put
\[
  {\mathcal U}=\bigcup_{i=1}^t{\mathcal U}_i,\qquad
  s=|{\mathcal T}|,\qquad u=|{\mathcal U}|.
\]
Their total size is $ts+u$.

First suppose ${\mathcal T}=\varnothing$. The families
${\mathcal U}_1,\ldots,{\mathcal U}_t$ are uncomplemented and strongly
incomparable, so Theorem~\ref{thm:hilton} gives $u\le2^{n-1}$.
Since $n\ge2$, we obtain $u\le2^{n-1}\le2^n-2\le tB_n$.

If ${\mathcal T}$ contains $\varnothing$ or $[n]$, then its antichain
property and strong incomparability with ${\mathcal U}$ imply $s=1$
and $u=0$, again giving the bound. We may therefore assume that
${\mathcal T}$ is non-empty and contains neither $\varnothing$ nor
$[n]$. It follows that ${\mathcal U}$ also contains neither of these
sets.

Write $\rho=L_n({\mathcal T})$. By
Corollary~\ref{lem:uncomplemented-lym},
\begin{equation}\label{eq:rho-s}
  \frac{s}{B_n}\le\rho\le1.
\end{equation}
For $1\le j\le n-1$, put
\[
  u_j=|{\mathcal U}\cap\binom{[n]}j|.
\]
For every $j$, the family
\[
  {\mathcal F}_j
  ={\mathcal T}\cup\left({\mathcal U}\cap\binom{[n]}j\right)
\]
is an antichain. If $n$ is odd, the usual LYM inequality gives
$\rho+u_j/\binom{n}{j}\le1$. If $n=2k$, every middle-level member of
${\mathcal T}$ belongs to $P({\mathcal F}_j)$, because
${\mathcal T}$ is uncomplemented and
$\overline{\mathcal T}\cap{\mathcal U}=\varnothing$.
Lemma~\ref{lem:unpaired-lym} therefore gives the same inequality:
\[
  \rho+\frac{u_j}{\binom{n}{j}}
  \le \sum_{F\in{\mathcal F}_j}\frac1{\binom{n}{|F|}}
       +\frac{|P({\mathcal F}_j)|}{k\binom{n}k}
  \le1.
\]
In particular, complementary pairs inside the middle level of
${\mathcal U}$ do not require discarding any sets. For both parities,
\begin{equation}\label{eq:noncentral-level}
  u_j\le\binom{n}{j}(1-\rho)
  \qquad(1\le j\le n-1).
\end{equation}
Summing over the non-empty proper levels gives
\begin{equation}\label{eq:unified-u-bound}
  u\le(2^n-2)(1-\rho)
    \le(2^n-2)\left(1-\frac{s}{B_n}\right).
\end{equation}
Consequently,
\begin{align*}
  ts+u
  &\le ts+(2^n-2)\left(1-\frac{s}{B_n}\right)\\
  &=tB_n+(2^n-2-tB_n)\left(1-\frac{s}{B_n}\right)\\
  &\le tB_n,
\end{align*}
since $2^n-2\le tB_n$ and $s\le B_n$.

Finally, the level $\binom{[n]}{\lfloor n/2\rfloor+1}$ is an
uncomplemented antichain of size $B_n$. Taking all $t$ families equal
to this level attains the bound.
\end{proof}

\section{Maximum size of non-empty, uncomplemented and strongly incomparable families}\label{sec:concluding}

 If empty families are allowed, Theorem~\ref{thm:hilton} is sharp by taking ${\mathcal A}_1$ and ${\mathcal A}_2$ as in Construction~\ref{cons-h} and taking all remaining families to be empty. It is interesting to ask what happens if every family must be non-empty.

\begin{problem}\label{prob:nonempty}
What is the maximum value of $\sum_{i=1}^t |{\mathcal A}_i|$ among all non-empty, uncomplemented and strongly incomparable families ${\mathcal A}_1,\dots,{\mathcal A}_t\subseteq 2^{[n]}$?
\end{problem}

For each pair $(n,t)$ for which such a system exists,  denote this maximum number by $M(n,t)$. It is clear that $M(n,1)=2^{n-1}$ for every $n\ge1$. Hilton's theorem and Construction~\ref{cons-h} imply that
 $M(n,2)=2^{n-1}$ for every $n\ge2$. For $n\ge3$ and $3\le t\le n$, we show that $M(n,t)$ is either $2^{n-1}-1$ or $2^{n-1}$.

\begin{proposition}\label{prop:nonempty-basic}
If $n\ge3$ and $3\le t\le n$, then
\[
  2^{n-1}-1\le M(n,t)\le2^{n-1}.
\]
\end{proposition}

\begin{proof}
The upper bounds follow from Theorem~\ref{thm:hilton}.  Set
\[
  {\mathcal A}_1=\{X\subseteq[n]:1\in X,\ |X|\le n-2\}.
\]
Partition the $(n-2)$-nd level of $2^{[n]\setminus\{1\}}$ into $t-1$ non-empty families:
\[
  {\mathcal A}_2\mathbin{\dot\cup}\cdots\mathbin{\dot\cup}{\mathcal A}_t
  =\binom{[n]\setminus\{1\}}{n-2}.
\]
Such a partition exists because this level has $n-1$ sets. The families are non-empty and uncomplemented. The last $t-1$ families lie in one level, so sets from different such families are incomparable. If $A\in{\mathcal A}_1$ and $B\in{\mathcal A}_i$ with $i\ge2$, then $1\in A\setminus B$, so $A\nsubseteq B$, while $|B|=n-2\ge|A|$ and $A\ne B$, so $B\nsubseteq A$. Thus the families are strongly incomparable.

Moreover,
\[
  |{\mathcal A}_1|=2^{n-1}-n,
  \qquad
  \sum_{i=2}^t|{\mathcal A}_i|=n-1.
\]
Hence their total size is $2^{n-1}-1$, which gives the lower bound.
\end{proof}

The proposition leaves a gap of one in the range $n\ge3$ and $3\le t\le n$. We conjecture the value $2^{n-1}-1$ is the maximum value.

\begin{conjecture}\label{conj:nonempty}
If $n\ge3$ and $3\le t\le n$, then
\[
  M(n,t)=2^{n-1}-1.
\]
\end{conjecture}

\appendix
\section{Appendix: The error in Theorem 1.5 of Goldwasser--Hilton--Zheng}\label{app:ghz-error}

Goldwasser, Hilton and Zheng \cite{GHZ19} stated that their weak-incomparability conjecture is equivalent to Hilton's strong-incomparability conjecture; this is their Theorem~1.5. The error occurs in the proof of the implication from the strong version to the weak version, more precisely in the proof of their Lemma~2.1.

We use their notation. Let $\mathcal S$ be the family of sets that occur in more than one of $\mathcal A_1,\ldots,\mathcal A_t$, let $\mathcal S_i=\mathcal A_i\cap\mathcal S$, and let
\[
  \mathcal T=(\mathcal A_1\cup\cdots\cup\mathcal A_t)\setminus\mathcal S.
\]
Their proof correctly observes that $\mathcal S_i$ and $\mathcal T$ are incomparable. It then claims, using their Lemma~1.17, that for a uniformly random maximal chain $\mathcal C$,
\begin{equation}\label{eq:ghz-false-step}
 \Pr\bigl(\mathcal C\cap(\mathcal S_i\cup\mathcal T)\ne\varnothing\bigr)
 \ge
 \frac{|\mathcal S_i|}{B_n}+\frac{|\mathcal T|}{2^{n-1}},
 \qquad
 B_n=\binom{n}{\lfloor n/2\rfloor+1}.
\end{equation}
This step is not valid. The chain-hitting lemma quoted as their Lemma~1.17 gives, for an arbitrary family $\mathcal F\subseteq2^{[n]}$, the lower bound
\[
  \Pr(\mathcal C\cap\mathcal F\ne\varnothing)\ge \frac{|\mathcal F|}{2^n}.
\]
It does not allow $2^n$ to be replaced by any smaller extremal upper bound for $|\mathcal F|$. Likewise, the fact that the uncomplemented antichain $\mathcal S_i$ satisfies $|\mathcal S_i|\le B_n$ does not imply that its chain-hitting probability is at least $|\mathcal S_i|/B_n$. Thus the improved cardinality bounds $|\mathcal S_i|\le B_n$ and $|\mathcal T|\le2^{n-1}$ do not give the two probability bounds used in \eqref{eq:ghz-false-step}. An extremal bound for the size of a family cannot simply be used as the denominator in a chain-hitting estimate.

Construction~A gives a direct way to see the problem. Take $h=1$ and write
\[
  W=[n]\setminus\{1\}.
\]
Then
\[
  \mathcal A_1=\{X\subseteq[n]:1\in X,\ |X|\le n-1\},
  \qquad
  \mathcal A_2=\cdots=\mathcal A_t=\{W\}.
\]
These families are uncomplemented and weakly incomparable. In the notation above,
\[
  \mathcal S=\{W\},\qquad
  \mathcal T=\mathcal A_1,
\]
and, for every $i\ge2$,
\[
  \mathcal S_i=\{W\},\qquad
  |\mathcal T|=2^{n-1}-1.
\]
Every maximal chain contains exactly one set of size $n-1$. If this set is $W$, then the chain meets $\mathcal S_i$; otherwise it contains $1$ and hence belongs to $\mathcal T$. Therefore
\[
  \Pr\bigl(\mathcal C\cap(\mathcal S_i\cup\mathcal T)\ne\varnothing\bigr)=1.
\]
On the other hand, the right-hand side claimed in \eqref{eq:ghz-false-step} is
\[
  \frac1{B_n}+\frac{2^{n-1}-1}{2^{n-1}}
  =1+\frac1{B_n}-\frac1{2^{n-1}}>1
  \qquad(n\ge3),
\]
because $B_n<2^{n-1}$. Thus \eqref{eq:ghz-false-step} cannot be true.

All later estimates in their Lemma~2.1, and hence their proof that Hilton's conjecture implies the weak-incomparability conjecture, depend on this step. Construction~A already shows that the weak-incomparability conjecture is false for every $t\ge3$, whereas Theorem~\ref{thm:hilton} proves Hilton's strong-incomparability conjecture. Therefore the equivalence asserted in Theorem~1.5 of \cite{GHZ19} is false.

\section*{Acknowledgement.}
We are grateful to Professor A.~J.~W.~Hilton for his lasting inspiration. In 1989, when Cheng Zhao was a graduate student under Professor Hilton's supervision, Professor Hilton mentioned this conjecture to him. That early discussion remained in memory and ultimately helped motivate the present work. 
This research was supported in part by the National Natural Science Foundation of China (Nos.~12601660 and 12571363) and by the Natural Science Foundation of Hunan Province (Grant No.~2025JJ30003).

\section*{Declaration on the use of generative AI}
 We used  ChatGPT5.6 to  find the errors in the proof  of Theorem 1.5 \cite{GHZ19} by Goldwasser, Hilton and Zheng(Section \ref{app:ghz-error} in the current version).
We  also used ChatGPT5.6 to  improve the language. We take full responsibility for the content of the paper.

\end{document}